\documentclass[12pt]{amsart}
\usepackage[all]{xy}
\usepackage{amsfonts}
\usepackage{amssymb}
\usepackage{pb-diagram}
\usepackage{amsmath}

\newtheorem{teo}{Theorem}[section]
\newtheorem{lem}[teo]{Lemma}
\newtheorem{prop}[teo]{Proposition}

\newtheorem{cor}[teo]{Corollary}

\newtheorem{dfn}[teo]{Definition}
\newtheorem{rk}[teo]{Remark}
\newtheorem{ex}[teo]{Example}

\begin{document}

\title[Module weak Banach-Saks and module Schur properties of Hilbert
$C^*$-modules]{Module weak Banach-Saks and module Schur properties of Hilbert
$C^*$-modules}

\author{Michael Frank}
\thanks{}
\address{Hochschule f\"ur Technik, Wirtschaft und Kultur (HTWK)
Leipzig, Fachbereich IMN, PF 301166, D-04251 Leipzig, F.R.
Germany}
\email{mfrank@imn.htwk-leipzig.de}
\urladdr{http://www.imn.htwk-leipzig.de/\~{}mfrank}

\author{Alexander Pavlov}
\thanks{Partially supported by the RFBR (grant 07-01-00046), by the
joint RFBR-DFG project (RFBR grant 07-01-91555 / DFG project ''K-Theory,
$C^*$-algebras, and Index theory''), and by the `Italian project Cofin06 -
Noncommutative Geometry, Quantum Groups and Applications'}
\address{Dipartimento di Matematica e Informatica, Universit{\`a}
degli Studi di Trieste, Piazzale Europa 1, I-34127 Trieste, Italy
and All-Russian Institute of Scientific and Technical Information,
Russian Academy of Sciences (VINITI RAS), Usievicha str. 20,
125190 Moscow A-190, Russia} \email{axpavlov@gmail.com}
\urladdr{http://www.axpavlov.com/}

\begin{abstract}
Continuing the research on the Banach-Saks and Schur properties
started in (cf. \cite{FP_BJMA}) we investigate analogous
properties in the module context. As an environment serves the
class of Hilbert $C^*$-modules. Some properties of weak module topologies
on Hilbert $C^*$-modules are described. Natural module analogues of the
classical weak Banach-Saks and Schur properties are defined and studied.
A number of useful characterizations of properties of Hilbert
$C^*$-modules is obtained. In particular, some interrelations of
these properties with the self-duality property
of countably generated Hilbert $C^*$-modules are established.
\end{abstract}

\maketitle
%%%%%%%%%%%%%%%%%%%%%%%%%%%%%%%%%
\newcounter{cou001}
%%%%%%%%%%%%%%%%%%%%%%%%%%%%%%%%%
% \tableofcontents
%%%%%%%%%%%%%%%%%%%%%%%%%%%%%%%%%
\section{Introduction}
%%%%%%%%%%%%%%%%%%%%%%%%%%%%%%%%%
In the present paper we extend the considerations of \cite{Kus:07,
FP_BJMA} on the Banach-Saks and Schur properties from the
situation of Banach spaces to the situation of certain Banach
$C^*$-modules. For Hilbert $C^*$-modules we introduce  some
natural module analogues of the weak Banach-Saks and Schur
properties and study their behavior. The main idea of these
constructions is to replace the weak topology of Banach spaces by
a weak topology of Hilbert $C^*$-modules. But there are at least
two suitable candidates for such a replacement - the weak topology
generated by the inner product and the weak topology generated by
$A$-linear bounded functionals on the Hilbert $C^*$-module. We
considered both of them, but mostly we are interested in the first
variant. In the third paragraph we obtain some properties of the
mentioned weak module topologies. The forth paragraph is dedicated
to the study of the module Schur property, where we investigate
such important cases as standard, countably generated and finitely
generated projective modules. In the fifth paragraph we focus on
the study of two variants of the weak module Banach-Saks property.
In the last paragraph some interrelations between both the module
Schur and the module weak Banach-Saks properties and self-duality
of Hilbert $C^*$-modules are considered.

%%%%%%%%%%%%%%%%%%%%%%%%%%%%%%%%%%%%%%%
\section{Preliminaries}
%%%%%%%%%%%%%%%%%%%%%%%%%%%%%%%%%%%%%%%
Let us recall that a (right) {\it pre-Hilbert $C^*$-module} over a
$C^*$-algebra $A$ (cf. \cite{Paschke}) is a right $A$-module $V$
equipped with an $A$-valued inner product $\langle \cdot , \cdot
\rangle : V \times V \to A$, which is $A$-linear in the second
variable, fulfils $ \langle x,y \rangle=\langle y,x
\rangle^{*}$; is {\it positive}, i.e. $\langle x,x \rangle \geq 0$,
and is {\it definite}, what means $\langle x,x \rangle = 0$ if and
only if $x=0$. A pre-Hilbert $A$-module is a {\it Hilbert $A$-module}
provided it is a Banach space with respect to the norm $\| x \|=
\|\langle x,x \rangle \|^{1/2}$. Hilbert modules over the field
of complex numbers $\mathbb{C}$ are Hilbert spaces.

Let $V_1$, $V_2$ be Hilbert $A$-modules. Then ${\rm
Hom}_A(V_1,V_2)$ stands for the set of all $A$-linear bounded
operators from $V_1$ to $V_2$. When $V_1=V_2=V$ we will write
${\rm End}_A(V)$ instead of ${\rm Hom}_A(V,V)$. An operator $T\in
{\rm Hom}_A(V_1,V_2)$ admits an adjoint operator $T^*\in {\rm
Hom}_A (V_2,V_1)$ if the equality
\begin{gather*}
    \langle Tx,y\rangle=\langle x, T^*y\rangle
\end{gather*}
holds for all $x\in V_1$, $y\in V_2$.  Let us remark that a module
operator between Hilbert $C^*$-modules may not have an adjoint one
unlike in the Hilbert space situation (cf.~\cite[Example
2.1.2]{MaTroBook}). In the sequel, by ${\rm End}^*_A(V)$ we will
denote the subset of ${\rm End}_A(V)$ which consists of operators
that possess an adjoint one. Both these Banach algebras can be
characterized alternatively in terms of different kinds of
multiplier algebras of the $C^*$-algebra ${\rm K}_A(V)$ of
''compact'' operators generated as the norm-closure of the linear
hull of the set of elementary operators $\{ \theta_{x,y} :
\theta_{x,y}(z) := x \langle y,z \rangle, \, x,y \in V \}$ on $V$,
cf.~\cite[Th.~15.2.12, 15.H]{W-O}.

To give some information on multiplier theory consider the
canonical isometric embedding of a $C^*$-algebra $A$ into its
bidual Banach space $A^{**}$, which may be equivalently described
as the enveloping von Neumann algebra of $A$ \cite{DixBook}.
Then (left, right) multipliers of $A$ may be defined as the
$C^*$-algebra $M(A)=\{a\in A^{**} : aA \subset A; Aa \subset A\}$
of \emph{multipliers} of $A$, and as the Banach algebra of
\emph{left} (resp., \emph{right}) \emph{multipliers}
as $LM(A)=\{a\in A^{**} : aA \subset A\}$, $RM(A)=LM(A)^*$.
More information about multipliers may be found in
\cite{Lance,MaTroBook,Pedersen,W-O}.

Resorting to the particular case of $C^*$-algebras $A$ which can
be considered as a right Hilbert $A$-module over themselves with
the inner product $\langle a,b\rangle=a^*b$ for $a,b \in A$, we
would like to introduce some topologies related to different kinds
of multipliers of $A$. The C*-algebra $M(A)$ of multipliers of $A$
can be identified with the set $M(A)={\rm End}^*_A(A)$, and the
Banach algebra $LM(A)$ of left multipliers of $A$ can be
alternatively defined as $LM(A)={\rm End}_A(A)$
(cf.~\cite{Kasparov, LinPJM}). These algebras may be understood
as maximal unitalizations of $A$ in some natural sense
(\cite{Lance, Pedersen, Pav1}), whereas the $C^*$-algebra $A_1$
with adjoint identity is a minimal unitalization of $A$.

For instance, in the commutative case, i.e. for $A=C_0(X)$, the
$C^*$-algebra of vanishing at infinity continuous functions on an
arbitrary locally compact Hausdorff space $X$, one has
$M(A)=LM(A)=C(\beta X)$ and  $A_1=C(X_+)$, where $\beta X$ and
$X_+$ are the Stone-\v{C}ech and the one-point compactifications
of $X$, respectively.

Keeping in mind that $A$ can be considered as
a $C^*$-subalgebra of ${\rm End}^*_A(A)$, we define the \emph{strict
topology} on $A$ (or even on ${\rm End}^*_A(A)$, if it is
desirable) as a locally convex topology induced by the family of
semi-norms $ \{\|\cdot\|_{a,l}:=\|\cdot a\|, \|\cdot\|_{a,r}:=\|a
\cdot\| : a\in A\}.$  Analogously, the \emph{left strict topology}
on $A$ (or even on ${\rm End}_A(A)$) is a locally-convex topology
induced by the family of semi-norms $\{\|\cdot\|_{a,l}:=\|\cdot
a\| : a\in A\}$. Analogously, introduce the \emph{right strict
topology} on $A$ as a locally convex topology induced by the family
of semi-norms $\{\|\cdot\|_{a,r}:=\|a \cdot\| : a\in A\}$.

The main property of these topologies is that the completion of
$A$ with respect to the strict (resp., left strict, right strict)
topology in $A^{**}$ gives the multiplier algebra $M(A)$ (resp.,
the left multiplier algebra $LM(A)$, the right multiplier algebra
$RM(A)$), (\cite{W-O, FP_JOP, Pedersen}).  The dual character of
$A^{**}$ as the universal $*$-representation of $A$ and as the
enveloping von Neumann algebra of $A$ makes sure, that the fact
remains true for any faithful non-degenerate $*$-representation of
$A$ in another $C^*$-algebra of operators in certain Hilbert
$C^*$-module replacing $A^{**}$ as an environment for completion
(cf. \cite{FP_JOP}).

For the convenience of the reader we also give a sketch of the
dual module theory that will be essential in the sequel.
Suppose $V$ stands for a Hilbert $A$-module again. Let
$V'={\rm Hom}_A(V,A)$ be the set of all $A$-linear
bounded maps from $V$ to $A$, the \emph{dual} (left Banach)
\emph{module} for $V$.
Apparently there is an isometric module embedding
\begin{gather}\label{eq:map_wedge}
   \wedge : V\rightarrow V',\quad x^{\wedge}(\cdot)=\langle
    x,\cdot\rangle.
\end{gather}
A Hilbert $A$-module is said to be \emph{self-dual} if this
map is surjective. Generally speaking, $V'$ is far not always a Hilbert
$C^*$-module provided $A$ is just an arbitrary $C^*$-algebra. The
inner product can be extended from $V$ to $V'$ (and even in such
way that the dual module will be self-dual) if and only if $A$ is
monotone complete (cf. \cite[Theorem 4.7]{FrankMN},
\cite{HamanaIJM, LinPJM}, although in the more important case of
$W^*$-algebras such an extension was described in \cite{Paschke}). As
usual, the \emph{bidual module} $V''$ should be define by the
second iteration, i.e. $V''={\rm Hom}_A(V',A)$. For a functional
$F\in V''$ one can define a functional
$\widetilde{F}\in V'$ by the formula
\begin{gather*}
    \widetilde{F}(x)=F(\widehat{x}),\quad x\in V \, .
\end{gather*}
Then the map $F\mapsto \widetilde{F}$ is an $A$-module isometry
from $V''$ to $V'$, and it allows to extend the inner product from
$V$  to $V''$ in the following way
\begin{gather*}
    \langle F,G\rangle=F(\widetilde{G}),\quad F,G\in V'' \, ,
\end{gather*}
where the norm arising from this inner product coincides
with the operator norm of $V''$. Now let $x\in V$, $f\in
V'$ and put
\begin{gather}\label{eq:map_dot}
     \dot{x}(f)=f(x)^*.
\end{gather}
Then the map (\ref{eq:map_dot}) is an isometric module embedding
of $V$ in $V''$ (see \cite{MaTroBook, Paschke2} for details).
It is interesting to observe that in opposition to the situation of
Banach spaces the third dual $V'''$ for $V$  is isomorphic to $V'$,
and the forth dual $V''''$ is isomorphic to $V''$, whereas the
Banach modules $V$, $V'$ and $V''$ may be pairwise non-isomorphic
in particular situations \cite{Paschke2}.

The \emph{standard Hilbert module} over a $C^*$-algebra $A$, that
is denoted by $l_2(A)$ or $H_A$, consists of all sequences $(a_i)$
with elements of $A$ such that the series $\sum_{i=1}^\infty
a_i^*a_i$ converge in norm. A Hilbert $C^*$-module $V$ is
\emph{countably generated} if there exists a sequence of elements
of $V$ such that the set of all finite $A$-linear combinations of
them is norm-dense in $V$. It is said to be \emph{finitely
generated} if for a finite subset of $V$ any element of $V$ can be
expressed as a finite $A$-linear combination of these generators.
Any finitely generated Hilbert $C^*$-module can be isometrically
embedded into a Hilbert $A$-module of type $A^n$ for some finite
$n \in \mathbb N$, \cite{MF}. Furthermore, Kasparov's
stabilization theorem asserts that any countably generated Hilbert
$A$-module can be represented as a direct orthogonal summand of
the standard module $l_2(A)$, \cite{Kasparov}. We need the
following facts about standard modules.

\begin{prop}\label{prop:dual for the standard module} {\rm (\cite{FrankZAA},
       \cite[Proposition 2.5.5]{MaTroBook})} \label{prop:funct_on_H_A}
  Consider the set of all sequences $f=(f_i)$, $f_i\in A$,
  $i\in \mathbb{N}$, such that the partial sums of the series
  $\sum f_i^*f_i$  are uniformly bounded. If $A$ is a unital
  $C^*$-algebra, then this set coincides with $H_A'$, where the
  action of $f$ on $H_A$ is defined by the formula
  \begin{displaymath}
     f(x)=\sum_{i=1}^\infty f_i^*x_i,
  \end{displaymath}
  where $x=(x_i)\in H_A$, and the norm of the functional $f$
  satisfies
  \begin{displaymath}
      \|f\|^2=\sup\limits_N\biggl\|\sum_{i=1}^\infty f_i^*f_i\biggr\|.
  \end{displaymath}
\end{prop}

\begin{prop}[\cite{FrankZAA, MaTroBook}]\label{prop:crit_selfd_H_A}
  Let $A$ be a $C^*$-algebra. Then the following conditions are
  equivalent:

  \begin{enumerate}%[\rm (i)]
    \item The Hilbert $C^*$-module $H_A$ is self-dual;
    \item The $C^*$-algebra $A$ is finite-dimensional.
  \end{enumerate}
\end{prop}

More information about Hilbert $C^*$-modules can be found
in \cite{Lance, MaTroJMS, MaTroBook,Raeburn-Williams,W-O},
for example.

%%%%%%%%%%%%%%%%%%%%%%%%%%%%%%%%%%%%%%%%%%%%%%%%
\section{Weak module topologies on Hilbert $C^*$-modules}
%%%%%%%%%%%%%%%%%%%%%%%%%%%%%%%%%%%%%%%%%%%%%%%%

We are going to establish some properties of the topologies
$\tau_{\widehat{V}}$ and $\tau_{V'}$ for arbitrary Hilbert
$C^*$-modules which are defined in the sequel. In particular, we
want to demonstrate that whenever both these topologies coincide
on a certain countably generated Hilbert module over a unital
$C^*$-algebra then this module has to be self-dual.

\begin{dfn}\rm Let $(V,\langle\cdot,\cdot\rangle)$ be a pre-Hilbert
  module over a $C^*$-algebra $A$, and let $V'$ be its dual module.
  The $V'$-weak module topology  $\tau_{V'}$ on $V$ is generated
  by the family of semi-norms
  \begin{gather}\label{eq:sem_top_dual}
    \{\nu_f\}_{f\in V'},\quad\text{where}\quad
    \nu_f(x)=\|f(x)\|,\, x\in V \, ,
  \end{gather}
  and the $\widehat{V}$-weak module topology $\tau_{\widehat{V}}$ on
  $V$ is generated by the family of semi-norms
  \begin{gather}\label{eq:ls_top_alg}
    \{\mu_z\}_{z\in V},\quad\text{where}\quad
    \mu_z(x)=\|\langle z,x\rangle\|,\, x\in V.
  \end{gather}
\end{dfn}

Obviously, $\tau_{\widehat{V}}$ is not stronger than $\tau_{V'}$,
and these topologies coincide whenever $V$ is self-dual. What about
the opposite interrelation? Do they always coincide?
The next example demonstrates that for arbitrary $V$ the
topologies $\tau_{\widehat{V}}$ and $\tau_{V'}$ are distinct,
in general.

\begin{ex}\rm
Let $(X,\rho)$ be a locally compact metric space, $A=C_0(X)$ be
the corresponding $C^*$-algebra and $V=A$ be the Hilbert
$A$-module under consideration. Suppose $\widetilde{X}=X\cup
\{\widetilde{x}\}$ stands for the one-point compactification of
$X$. By the commutativity of $A$ the multiplier algebras $M(A)$
and $LM(A)$ coincide. Moreover, the $A$-dual Banach $A$-module
$V'$ of $V$ can be identified with the Banach space ${\rm
End}_A(A)$ which is isometrically isomorphic to the left
multiplier algebra $LM(A)$ of $A$. So $V'$ exactly coincides with
the algebra $C_b(X)$ of all bounded continuous functions on $X$,
and $V'=C_b(X)$ can be considered to act by left multiplication on
$V=C_0(X)$ (see \cite{Murphy, W-O}). Consider a sequence $\{x_i\}$
in $X$ which converges to $\widetilde{x}$ with respect to the
topology of $\widetilde{X}$. For this sequence choose a sequence
of open neighborhoods
\[
   O_{\delta_i}(x_i)=\{y\in X : \rho (y,x_i)<\delta_i\} \, ,
\]
under the additional request that the sequence of positive numbers
$\{\delta_i\}$ has to be a null sequence. Let $f_i\in C_0(X)$ be
elements such that $f_i(x_i)=1$ and $f_i=0$ outside of
$O_{\delta_i}(x_i)$. Then for any $g\in C_0(X)$ one has
\begin{eqnarray*}
   \lim_{i \to \infty} \|gf_i\|
      & = &  \lim_{i \to \infty} \sup_{x\in X}|g(x)f_i(x)|
        =  \lim_{i \to \infty} \sup_{x\in O_{\delta_i}(x_i)}|g(x)f_i(x)|\\
     &\le & \lim_{i \to \infty} \sup_{x\in    O_{\delta_i}(x_i)}|g(x)|
     \, =  \, 0
\end{eqnarray*}
what means $\{f_i\}$ is a null sequence with respect to the
$\tau_{\widehat{V}}$-topology. On the other hand
consider the constant function $h$ on $X$, which always equals
to 1. Then $h$ belongs to $C_b(X)=V'$ and
\[
   \|hf_i\|=\sup_{x\in X}|f_i(x)|=1
\]
for any $i \in {\mathbb N}$. Consequently, $\{f_i\}$ does not
converge to zero with respect to the topology $\tau_{V'}$.
\end{ex}

\begin{prop}\label{prop:conv_bound1}
   Let $\{x_n\}$ be a $\tau_{\widehat{V}}$-converging
   sequence of a Hilbert $C^*$-module $V$. Then the set $\{x_n\}$
   is bounded in norm.
\end{prop}

\begin{proof}
Consider the sets
\[
   A_{k,n}=\{y\in V : \|\langle y,x_n\rangle\|\le k\}
\]
and $A_k=\cap_n A_{k,n}$. These sets are closed in norm, because
the function $y\mapsto \|\langle y,x_n\rangle\|$ is continuous for
any fixed $x_n$. The $\tau_{\widehat{V}}$-convergence of $\{x_n\}$
 implies that $V=\cup_k A_k$. By the Baire's theorem
there are a number $k_0$ and a ball $B(y,\varepsilon)=\{z\in V :
\|y-z\|<\varepsilon\}$ such that $B(y,\varepsilon)\subset
A_{k_0}$. Thus the sequence $\{x_n\}$, understood as a subset
of $V'$, is bounded on the ball $B(y,\varepsilon)$ and,
consequently, on any ball of $V$, in particular on $B(0,1)$, that
exactly means
\begin{eqnarray*}
    C & \ge & \sup_{\|y\|\le 1} \|\langle y,x_n\rangle\|
         =    \sup_{\|y\|\le 1} \|\widehat{y}(x_n)\|\\
      &  =  & \sup_{\|y\|\le 1} \|\dot{x}_n(\widehat{y})\|
         =    \sup_{\|\widehat{y}\|\le 1} \|\dot{x}_n(\widehat{y})\| \\
      &  =  & \|\dot{x}_n\|
         =    \|x_n\|
\end{eqnarray*}
for some constant $C$ and for all $n \in {\mathbb N}$, where the properties
of the maps (\ref{eq:map_wedge}) and (\ref{eq:map_dot}) have been applied.
\end{proof}

\begin{prop}
   Let $\{x_n\}$ be a $\tau_{V'}$-converging sequence of a pre-Hilbert
   $C^*$-module $V$. Then the set $\{x_n\}$ is bounded in norm.
\end{prop}

\begin{proof} We should only remark that $V'$ is a Banach module even
if $V$ is just a pre-Hilbert one. Further arguments are close
to those applied in the proof of Proposition \ref{prop:conv_bound1}.
\end{proof}

Let us remind that a subset $M$ of a topological vector space $X$ is
said to be \emph{bounded} if for any open neighborhood $U$ of zero in
$X$ there is real number $s>0$ such that $M\subset tU$ for any $t>s$.

\begin{prop}\label{prop:conv_bound3}
   The following conditions for a subset $Q$ of a Hilbert
   $C^*$-module $V$ are equivalent:
   \begin{enumerate}
    \item $Q$ is bounded with respect to the topology $\tau_{\widehat{V}}$;
    \item $Q$ is norm-bounded.
   \end{enumerate}
\end{prop}

\begin{proof} Because the norm topology is stronger than the topology
$\tau_{\widehat{V}}$ one has only to check the implication ${\mathrm{(i)}}
\rightarrow {\mathrm{(ii)}}$.

Let $Q$ be bounded with respect to the topology
$\tau_{\widehat{V}}$. Assume there is a sequence $\{x_n\} \in Q$
for which the sequence $\{ \|x_n \| \}$ is unbounded. Because the
set $\{x_n\}$ is $\tau_{\widehat{V}}$-bounded, for any $y\in V$
one can find a real number $t>0$ satisfying
\begin{gather*}
    \{x_n\}\subset t\{x : \|\langle y,x\rangle\|<1\}
\end{gather*}
what yields
\begin{gather*}
   \|\langle y,x_n\rangle\|<t\quad\text{for all}\quad n\in
   \mathbb{N}.
\end{gather*}
In particular, the sets $\{\|\langle y,x_n\rangle\|\}_{n=1}^\infty$ are
bounded for any fixed $y\in V$. The remainder of the proof can be done
in the same way like the proof of Proposition \ref{prop:conv_bound1}.
We arrive at a contradiction since the sequence $\{ \|x_n \| \}$
turns out to be bounded.
\end{proof}

\begin{prop}
  The following conditions for a subset $Q$ of a pre-Hilbert $C^*$-module
  $V$ are equivalent:
  \begin{enumerate}
    \item $Q$ is bounded with respect to the topology $\tau_{V'}$;
    \item $Q$ is norm-bounded.
  \end{enumerate}
\end{prop}

\begin{proof}
The proof follows the arguments of the proofs of Propositions
\ref{prop:conv_bound1} and \ref{prop:conv_bound3} with obvious
changes in details.
\end{proof}

\begin{teo}\label{teo:coincidence_topol}
  Let $A$ be a unital $C^*$-algebra, $V=l_2(A)$ be the
  standard countably ge\-nerated Hilbert $A$-module.
  Suppose, the topologies $\tau_{\widehat{V}}$ and $\tau_{V'}$ coincide
  on $V$. Then $V$ is self-dual and $A$ has to be finite-dimensional.
\end{teo}

\begin{proof} Let us take into consideration any functional $\beta\in V'$.
Then by Proposition \ref{prop:crit_selfd_H_A} there are elements
$b_i$ from $A$ such that $\beta (x)=\sum_{i=1}^\infty b_i^* a_i$,
where $x=(a_i)\in V$ and
\[
\biggl\|\sum_{i=1}^N b_i^*b_i\biggr\|\le C
\]
for some constant $C$ and for all $N$. By the suppositions of the
theorem there are vectors $y_1,\dots, y_n$ of $V$ and a constant
$K>0$ such that
\begin{gather}\label{eq:ineq_nu_mu}
    \nu_{\beta}(x)\le K\max\{\mu_{y_1}(x),\dots,\mu_{y_n}(x)\}
\end{gather}
for any $x\in V$, where we have used the notations
(\ref{eq:sem_top_dual}) and (\ref{eq:ls_top_alg}). Let us fix the
notations $y_i=(y^{(i)}_j)$, $x=(a_i)$, where $y^{(i)}_j, a_i\in
A$ and rewrite the inequality (\ref{eq:ineq_nu_mu}) in the form
\begin{gather}\label{eq:ineq_nu_mu_2}
    \biggl\|\sum_{i=1}^\infty b_i^*a_i\biggr\|\le
    K\max\left\{\biggl\|\sum_{j=1}^\infty \left(y_j^{(1)}\right)^*a_j\biggr\|,\dots,
    \biggl\|\sum_{j=1}^\infty \left(y_j^{(n)}\right)^*a_j\biggr\|\right\}.
\end{gather}
For arbitrary positive integers $k$ and $N$ set
$a_i=b_i$ if $k\le i\le k+N$, and $a_i=0$ otherwise. Then the
formula (\ref{eq:ineq_nu_mu_2}) gives
\begin{align*}
  \biggl\|\sum_{i=k}^{k+N} b_i^*b_i\biggr\|\le K\max
  &\left\{\biggl\|\sum_{j=k}^{k+N}
  \left(y_j^{(1)}\right)^*b_j\biggr\|,\dots,
    \biggl\|\sum_{j=k}^{k+N}
    \left(y_j^{(n)}\right)^*b_j\biggr\|\right\}\\
    \le  K\max &\left\{\biggl\|\sum_{j=k}^{k+N}
  \left(y_j^{(1)}\right)^*y_j^{(1)}\biggr\|^{1/2}
  \biggl\|\sum_{j=k}^{k+N} b_j^*b_j\biggr\|^{1/2},\dots,\right.\\
   &\left.   \biggl\|\sum_{j=k}^{k+N}
  \left(y_j^{(n)}\right)^*y_j^{(n)}\biggr\|^{1/2}
  \biggl\|\sum_{j=k}^{k+N} b_j^*b_j\biggr\|^{1/2}\right\}\\
  \le C^{\frac{1}{2}} \max &\left\{\biggl\|\sum_{j=k}^{k+N}
  \left(y_j^{(1)}\right)^*y_j^{(1)}\biggr\|^{1/2} ,\dots,
  \biggl\|\sum_{j=k}^{k+N}
  \left(y_j^{(n)}\right)^*y_j^{(n)}\biggr\|^{1/2} \right\}.
\end{align*}
Therefore $\{\sum_{i=1}^{N} b_i^*b_i\}_{N=1}^\infty$ is a Cauchy
sequence with respect to the uniform topology what ensures the
vector $(b_i)$ to belong to $V=l_2(A)$. Thus, $V$ is self-dual
and, by Proposition \ref{prop:crit_selfd_H_A}, $A$ is
finite-dimensional.
\end{proof}

\begin{cor}
  Let $A$ be a unital $C^*$-algebra. Let $V$ be a countably
  generated Hilbert $A$-module, for which the topologies
  $\tau_{\widehat{V}}$ and $\tau_{V'}$ coincide. Then $V$ is
  self-dual.
\end{cor}

\begin{proof} By Kasparov's stabilization theorem
\cite{Kasparov} $V$ can be isometrically represented as a direct
orthogonal summand of $l_2(A)$, also denoted by $V$ for brevity.
So the Hilbert $A$-module $l_2(A)$ decomposes into a direct
orthogonal sum $V \oplus V^\bot$. This orthogonal decomposition of
$l_2(A)$ causes an direct topological sum decomposition of the
$A$-dual Banach $A$-module $l_2(A)'$ into the sum of the set
$(V^\bot)'$ of all elements of $l_2(A)'$ vanishing on $V$ and the
set  $V'$ of all elements of $l_2(A)'$ vanishing on $V^\bot$.
Applying the basic idea of the previous proof again and resorting
to sequences in $V \subseteq l_2(A)$ the standard isometric
embedding of $l_2(A)$ into its $A$-dual $l_2(A)'$ is seen to map
$V$ onto $V'$. Consequently, $V$ is self-dual.
\end{proof}

Summing up, we have shown that for countably generated
Hilbert $C^*$-modules over unital $C^*$-algebras the coincidence
of both the toplogies $\tau_{\widehat{V}}$ and $\tau_{V'}$
on the Hilbert $C^*$-module forces self-duality of the
Hilbert $C^*$-module, and vice versa.

%%%%%%%%%%%%%%%%%%%%%%%%%%%%%%%%%%%%
\section{The module Schur property}
%%%%%%%%%%%%%%%%%%%%%%%%%%%%%%%%%%%

In the following section we introduce a module analogue
of the Schur property of Banach spaces. Recall that a Banach
space $X$ has the \emph{Schur property} if every weak convergent
sequence in $X$ converges in norm. Let $Y$ be a closed subspace
of $X$. Then $X$ has the Schur property if and only if both $Y$
and the quotient space $X/Y$ have the same property (cf.
\cite{CasGon, Kus:07}).

\begin{ex}{\rm
The Banach space $L^1[0,1]$ has the Schur property. Indeed, let
$\{x_i\}$ be a weakly convergent sequence of $L^1[0,1]$. Then, in
particular, for a unit function $1\in L^\infty[0,1]=L^1[0,1]'$ the
sequence $\{1\cdot x_i\}$ has to converge in norm, but, in fact,
this sequence coincides with $\{x_i\}$. }
\end{ex}

Comparing the Banach space and the Banach module variants of the
Schur properties we establish properties of $C^*$-algebras and of
Hilbert $C^*$-modules with the module Schur property. For
classification results with respect to the classical Schur
property see \cite{Kus:07,FP_BJMA}.

\begin{dfn}\label{def:Schur_prop}\rm
  A Hilbert $C^*$-module $V$ has the \emph{module Schur
  property} if every $\tau_{\widehat{V}}$-convergent
  sequence in $V$  converges in norm.
\end{dfn}

\begin{lem}\label{lem:seq_compl_Schur}
  Any $\tau_{\widehat{V}}$-convergent sequence  $\{ x_n \}$
  of a Hilbert $C^*$-module $V$with a subsequence $\{ x_{n(l)} \}$
  which converges in norm to some element $x \in V$ admits $x$
  as its $\tau_{\widehat{V}}$-limit.  Consequently, Hilbert
  $C^*$-modules with the module Schur property are sequentially
  $\tau_{\widehat{V}}$-complete.
  \end{lem}

\begin{proof}
For any $z \in V$ and any $\varepsilon > 0$ there exists a number
$N \in \mathbb N$ such that
\[
   \| \langle z, x_n-x_m \rangle \| <  \varepsilon
\]
for any $m,n > N$. In particular, there exists a number $L \in
\mathbb N$ such that $n(l) > N$ for any $l > L$. So
\begin{equation*}
   \| \langle z, x-x_m \rangle \|  =
   \lim_{l \to \infty} \| \langle z, x_{n(l)}-x_m \rangle \|
        \, \leq  \, \varepsilon
\end{equation*}
for any $m > N$, and $x \in V$ turns out to be the
$\tau_{\widehat{V}}$-limit of the sequence $\{ x_n \}$.
\end{proof}

Let us recall that a $C^*$-algebra has the Schur property if and
only if it is finite dimensional (cf.~\cite[Lemma 3.8]{Kus:07}).
Our next goal is to prove that a $C^*$-algebra with a strictly
positive element possesses the module Schur property if and only
if it is unital. The conclusion that a unital $C^*$-algebra has
the module Schur property is obvious. But, before investigating
the converse statement we will consider a couple of examples as a
motivation.

\begin{ex}\rm Let $A=V=C_0(0,1]$ and let us define a function $f_n
\in A$ in the following way: $f$ equals to zero at $0$ and on the
interval $[\frac{1}{n},1]$, $f$ equals to $1$ at the point $\frac{1}{2n}$,
and $f$ is linear on both the intervals $[0,\frac{1}{2n}]$ and
$[\frac{1}{2n},\frac{1}{n}]$. Obviously, the sequence
$\{f_n\}$ converges to zero with respect to the topology
$\tau_{\widehat{V}}$, but it does not converge in norm. So the
$C^*$-algebra $C_0(0,1]$ does not have the module Schur property.
\end{ex}

\begin{ex}\rm \label{ex:comp_op}
Let $A=V=K(H)$ be the $C^*$-algebra  of  compact operators in a
separable Hilbert space $H$ and $B(H)$ be the set of all bounded
linear operators on $H$. Consider a sequence $\{a_n\} \in A$ and
assume that the sequence $\{ka_n\}$ converges in norm for any
$k\in K(H)$. Then its limit with respect to the right strict
topology belongs to the closure of $K(H)$ in $B(H)$, and this
closure coincides with $B(H)$ (see, for instance, \cite{LinJOP,
Pav1}). However, the norm limit of $\{a_n\}$, whenever it exists,
has to belong to $K(H)$. So any sequence from $K(H)$, whose right
strict limit lies outside $K(H)$, $\tau_{\widehat{V}}$-converges,
but does not converge in norm. Therefore $K(H)$ does not have the
module Schur property, too.
\end{ex}

\begin{teo}\label{teo:Schur for C*-alg}
  A $C^*$-algebra with a strictly positive element, in particular,
  a separable   $C^*$-algebra (cf. \cite[1.4.3]{Pedersen}), has the module
  Schur property if and only if it is unital. In particular, the
  Schur property for $C^*$-algebras is strictly stronger than the
  module Schur property.
\end{teo}

\begin{proof} Suppose a non-unital $C^*$-algebra $A$ has a strictly
positive element. That condition on $A$ is equivalent to the
existence of a countable approximative unit in $A$
(\cite[Proposition 3.10.5]{Pedersen}).
For such $A$ any element of its right multiplier algebra $RM(A)$
may be obtained just as a limit of a sequence of $A$ converging
with respect to the right strict topology, i.e. there is not
any necessity to consider nets under these
assumptions (cf. \cite[\S 2.3]{W-O}, \cite[\S 5.5]{MaTroBook}).
Now, to check that $A$ does not have the module Schur property
we will reason similar as at Example~\ref{ex:comp_op}. Let a
sequence $\{a_n\}$ of $A$ be $\tau_{\widehat{V}}$-convergent
(here we identify $V=A$ again), i.e. the sequence $\{b^*a_n\}$
converges in norm for any $b \in A$. Therefore, the
$\tau_{\widehat{V}}$-limit of $\{a_n\}$ belongs to the right strict
closure  of $A$ inside of its bidual Banach space $A^{**}$,
i.e. it belongs to the right multipliers $RM(A)$ of $A$.
Because the algebra $RM(A)$ is strictly larger than $A$ whenever
$A$ is not unital, we can find a sequence from $A$ that converges
with respect to the topology $\tau_{\widehat{V}}$ (for instance,
to the identity of $RM(A)$), but does not converge in norm. Thus $A$
does not have the module Schur property.

Since unital $C^*$-algebras have the Schur property if and only if
they are finite-dimensional, the module Schur property is
weaker than the Schur property for $C^*$-algebras.
\end{proof}

%Because any separable $C^*$-algebra has a countable approximative
%unit (cf. \cite[1.4.3]{Pedersen}) the next assertion holds.

%\begin{cor}
%  Any separable non-unital $C^*$-algebra has the module Schur
%  property.
%\end{cor}

The next example demonstrates the requirement of Theorem
\ref{teo:Schur for C*-alg} to a $C^*$-algebra to have a strictly
positive element is essential.

\begin{ex}{\rm
Let $H$ be a non-separable Hilbert space, $A=B(H)$ be
the $C^*$-algebra of all linear bounded operators in $H$ and $A_0$
denotes the set of operators $T\in B(H)$ such that both
$\mathrm{Ker} (T)^\bot$ and $\overline{\mathrm{Range} (T)}$ are
separable (closed) subspaces of $H$. Then $A_0$
is an involutive subalgebra of $A$, containing all compact
operators on $H$. Moreover, $A_0$ is closed in norm. Indeed, let a
sequence $\{T_n\}$ of $A_0$ converges in norm to a certain element
$T$ of $B(H)$. Consider (closed) separable subspaces
\begin{gather}\label{eq:example_nonseparable_space}
   H_1=\overline{\mathrm{span}}\{\cup \mathrm{Ker} (T_n)^\bot\},\quad
   H_2=\overline{\mathrm{span}}\{ \cup \overline{\mathrm{Range} (T_n)}\}
\end{gather}
 of $H$, generated by the unions $\cup
\mathrm{Ker} (T_n)^\bot$ and $\cup \overline{\mathrm{Range}
(T_n)}$ respectively. Then, obviously, $H_1^\bot$ coincides with
the intersection $\cap \mathrm{Ker} (T_n)$, what infers $H_1^\bot$
belongs to $\mathrm{Ker} (T)$. This implies that
 $\mathrm{Ker} (T)^\bot$ is a subspace of $H_1$ and, consequently, is
 separable. On the other hand,  $\overline{\mathrm{Range} (T)}$
belongs to the union $\cup \overline{\mathrm{Range} (T_n)}$.
Therefore, both $\mathrm{Ker} (T)^\bot$ and
$\overline{\mathrm{Range} (T)}$ are separable subspaces of $H$ and
$T\in A_0$. Clearly, the $C^*$-algebra $A_0$ is not unital, but we
claim that it is not even $\sigma$-unital. To verify this,
let us suppose the opposite assertion to be true, and there is a
countable approximative identity $\{T_n\}$ in $A_0$. Let $H_1$ be
defined by (\ref{eq:example_nonseparable_space}). Then there is a
non-zero vector $x\in H$  such that it is orthogonal to $H_1$. Let
$p\in B(H)$ be the orthogonal projection onto the span of $x$.
Then, actually, $p$ belongs to $A_0$ and the sequence $pT_n$ does
not converge in norm to $p$. So we have a contradiction and,
consequently, $A_0$ does not have any countable approximative
identity. Nevertheless, we claim that $V=A_0$, considered as a
Hilbert $A_0$-module, satisfies the module Schur property. Indeed,
assume $\{T_n\}$ is a $\tau_{\widehat{V}}$-convergent sequence of
$V$. Then both spaces $H_1$ and $H_2$, defined by
(\ref{eq:example_nonseparable_space}), are separable. Let $u\in
B(H)$ be a minimal partial isometry between $H_1$ and $H_2$. Then, in
fact, $u$ belongs to $A_0$, where $u|_{H_1^\bot}=0$ and
$u^*|_{H_2^\bot}=0$. Thus the sequence $\{u^*T_n\}$ converges in
norm and coincides element-wise with the sequence $\{T_n\}$.}
\end{ex}

\begin{prop}\label{prop:direct_sum_Schur}
  Assume a Hilbert $C^*$-module $V$ is represented as a direct sum
  $V=V_1\oplus V_2$. Then the following conditions are equivalent:
  \begin{enumerate}
    \item $V$ has the module Schur property;
    \item Both $V_1$ and $V_2$ have the module Schur property.
  \end{enumerate}
\end{prop}

\begin{proof}
Obviously, (i) implies (ii) and we have just to verify the inverse
assertion. So consider any $\tau_{\widehat{V}}$-convergent sequence
$\{x_i\}$ of $V$. Then $x_i=x_i^{(1)}\oplus x_i^{(2)}$ with
$x_i^{(1)}\in V_1$, $x_i^{(2)}\in V_2$. Obviously, the sequences
$\{x_i^{(j)}\}$ are $\tau_{\widehat{V_j}}$-convergent sequences
($j=1, 2$), and therefore, they converge with respect to the norm
of $V_j$ ($j=1, 2$). Consequently,
\[
  \|x_n-x_m\|\le
  \|x_n^{(1)}-x_m^{(1)}\|+\|x_n^{(2)}-x_m^{(2)}\|
\]
for any $n, m$. Thus, $\{x_i\}$ is a Cauchy sequence with respect
to the norm topology of $V$.
\end{proof}

\begin{cor}\label{cor:fin_gen_mod_Schur}
  Finitely generated Hilbert modules over unital $C^*$-algebras $A$
  have the module Schur property. In particular, the standard
  Hilbert $A$-modules $A^n$ for $n \in \mathbb N$ have the module
  Schur property.
\end{cor}

\begin{proof} For any Hilbert $A$-module $A^n$ and for any $n
\in \mathbb N$ we can apply Proposition \ref{prop:direct_sum_Schur}
$n$ times, since $A$ is shown to have the module Schur property
by Theorem \ref{teo:Schur for C*-alg}.

By the Kasparov stabilization theorem (cf. \cite{Kasparov},
\cite[Theorem 1.4.2]{MaTroBook}) any finitely generated module $V$
is an orthogonal direct summand of the standard module $l_2(A)$.
Therefore, by~\cite[Theorem 2.7.5]{MaTroBook} (or by \cite{MF}) it
has to be projective whenever $A$ is unital. So, $V$ is a direct
summand of some Hilbert $A$-module $A^n$, $n < \infty$, and
Proposition \ref{prop:direct_sum_Schur} applies in the inverse
direction. So any projective finitely generated Hilbert module
over a unital $C^*$-algebra has the module Schur property.
\end{proof}

\begin{cor}\label{cor:fin_gen_proj_mod_Schur}
  Let $A$ be a $\sigma$-unital $C^*$-algebra with the module Schur
  property. Then finitely generated projective Hilbert modules over
  $A$ have the module Schur property.
\end{cor}

\begin{proof}
The assertion directly follows from Theorem \ref{teo:Schur for C*-alg} and Proposition
\ref{prop:direct_sum_Schur}.
\end{proof}

\begin{prop}\label{prop:Schur stand mod}
  The standard Hilbert module $V=l_2(A)$ does not have the module
  Schur property for any $C^*$-algebra $A$.
\end{prop}

\begin{proof}
For a unital $C^*$-algebra the statement is clear, because the
standard basis $\{e_i\}$ (with all entries of $e_i$ equal to
zero except the $i$-th, which equals $1_A$)
$\tau_{\widehat{V}}$-converges to zero, but does not converge in
norm.

For an arbitrary $C^*$-algebra $A$ we can reason in the following
way: let us fix an element $a\in A$ of norm one and consider the
sequence $\{x_k\}$ from $l_2(A)$, where all entries of $x_k$ are
zero except the $k$-th entry that equals to $a$. Then for any
$y=(b_i)$ from $l_2(A)$ one has
\[
   \|\langle y,x_k\rangle\|=\|b_k^*a\|\le\|b_k\|
\]
for each $k$, so $\{x_k\}$ converges to zero with respect to the
topology $\tau_{\widehat{V}}$. On the other hand, this sequence
does not converge in norm.
\end{proof}

Let us consider one example of a countably, but not finitely generated
and non-standard Hilbert $C^*$-module without the module Schur property.

\begin{figure}[h]
  \begin{picture}(90,60)%(95,0)
  \put(0,0){\line(1,0){70}} \put(80,5){$X$}
  \put(0,40){\line(1,0){40}} \put(40,40){\line(3,2){30}}
  \put(40,40){\line(3,-2){30}}  \put(80,50){$Y$}
  \put(30,30){\vector(0,-1){15}} \put(20,20){$p$}
  \end{picture}
 \caption{Example~\ref{ex:count_gen_mod}}\label{fig:count_gen_mod}
\end{figure}

\begin{ex}\label{ex:count_gen_mod}\rm
Let us consider the projection map $p: Y\rightarrow X$ from Figure
\ref{fig:count_gen_mod}, where $X$ is an interval, say $[-1,1]$,
and $Y$ is the topological union of one interval with two copies
of another half-interval with a branch point at $0$.
Then $C(Y)$ is a Banach $C(X)$-module with respect to the action
\[
   (f\xi)(y)=f(y)\xi(p(y)), \qquad f\in C(Y), \xi\in C(X).
\]
Let us define the $C(X)$-valued inner product on
$C(Y)$ by the formula
\begin{gather}\label{eq:inner_pr}
    \langle f,g\rangle (x)=\frac{1}{\# p^{-1}(x)} \sum_{y\in
    p^{-1}(x)} \overline{f(y)}g(y),
\end{gather}
where $\# p^{-1}(x)$ is the cardinality of $p^{-1}(x)$. It was
shown in  \cite{PavTro_cov} that $C(Y)$ is a countably, but not
finitely generated Hilbert $C(X)$-module with respect to the inner
product (\ref{eq:inner_pr}). Now for a point $x\in [0,1]$ let us
denote by $y_x^{(1)}$ its pre-image $p^{-1}(x)$ intersected with
the upper line of $Y$, and by $y_x^{(2)}$ the intersection of
$p^{-1}(x)$ with the lower line of $Y$. Consider functions $h_n\in
C(X)$ that equal zero to the left of the point $\frac{1}{2n}$ and
to the right of the point $\frac{1}{n}$, equal $1$ at the point
$\frac{3}{4n}$ and linear on both intervals
$[\frac{1}{2n},\frac{3}{4n}]$ and $[\frac{3}{4n},\frac{1}{n}]$.
Now we define a sequence $\{f_n\}$ from $C(Y)$ in the following
way: $f_n=0$ on $p^{-1}([-1,0])$, $f_n(y_x^{(1)})=h_n(x)$ and
$f_n(y_x^{(2)})=-h_n(x)$ for all $x\in [0,1]$. Then for any $g\in
C(Y)$ one has
\begin{align*}
    \|\langle g,f_n\rangle\|&=\max_{x\in
    [\frac{1}{2n},\frac{1}{n}]} \frac{1}{2}\left|
    \overline{g(y_x^{(1)})}f_n(y_x^{(1)})+
    \overline{g(y_x^{(2)})}f_n(y_x^{(2)})\right|\\
    &=\max_{x\in
    [\frac{1}{2n},\frac{1}{n}]}\frac{1}{2}\left|{g(y_x^{(1)})}-
    {g(y_x^{(2)})}\right|h_n(x)\\&\le \max_{x\in
    [\frac{1}{2n},\frac{1}{n}]}\frac{1}{2}\left|{g(y_x^{(1)})}-
    {g(y_x^{(2)})}\right|
\end{align*}
and the latter sequence converges to zero if $n$ goes to infinity.
But on the other hand the sequence $\{f_n\}$ does not converge in
norm. Thus $C(Y)$ does not have the module Schur property.
\end{ex}

%%%%%%%%%%%%%%%%%%%%%%%%%%%%%%%%%%%%%%%%%%%%%%%%%%%%%%%%%%%
\section{Module weak Banach-Saks properties }
%%%%%%%%%%%%%%%%%%%%%%%%%%%%%%%%%%%%%%%%%%%%%%%%%%%%%%%%%%%

The basic idea of the investigations presented in this section
is a search for module analogues for the different kinds
of Banach-Saks properties for Banach spaces. As most promissing
we select a certain generalization of the weak Banach-Saks
property.

Let us start with a short review of two Banach-Saks type properties
for the classical situation. A Banach space $X$ has the
\emph{Banach-Saks property} if from any bounded sequence
$\{x_n\}$ of $X$ there may be extracted a subsequence
$\{x_{n(k)}\}$ such that
\begin{gather}\label{eq:Banach-Saks}
    \lim_{k\rightarrow\infty}\biggl\|\frac{1}{k}\sum_{i=1}^k
    x_{n(i)} - x \biggr\|=0
\end{gather}
for some element $x \in X$.

Originally, this condition was introduced and studied by S.~Banach
and S.~Saks in \cite{B-S:30} just for particular examples, namely
they demonstrated that the Banach spaces $L^p([0,1])$ with $1 < p
< \infty$ satisfy the Banach-Saks property. Let us emphasize also
that a Banach space is reflexive whenever it has
the Banach-Saks property \cite{Die:75}. A $C^*$-algebra has the
Banach-Saks property if and only if it is finite-dimensional
\cite{Kus:07,Chu:94}.

In Functional Analysis a weaker type of the Banach-Saks property
gives very useful results for larger classes of Banach spaces.
More precisely, a Banach space $X$ has the \emph{weak Banach-Saks
property} if for any sequence $\{x_n\}$ of $X$ which converges
weakly to zero there may be selected a subsequence $\{x_{n(k)}\}$
such that the equality (\ref{eq:Banach-Saks}) is fulfilled with
$x=0$. We are interested  in certain module analogues of the weak
Banach-Saks property in this section. To give definitions, assume
$V$ to be a Hilbert $C^*$-module. Consider the following
conditions on $V$.
\begin{itemize}
    \item[{(mBS1)}] \label{eq:mwBS1}
       Any $\tau_{\widehat{V}}$-null sequence $\{x_n\}$ of $V$
       admits a subsequence $\{x_{n(i)}\}$ that satisfies the
       equality~(\ref{eq:Banach-Saks}) with $x=0$.
\end{itemize}
\begin{itemize}
    \item[{(mBS2)}] \label{eq:mwBS2}
       Any $\tau_{\widehat{V}}$-convergent sequence $\{x_n\}$ of
       $V$ admits a subsequence $\{x_{n(i)}\}$ such that the
       sequence $\frac{1}{k}\sum_{i=1}^k x_{n(i)}$ converges in
       norm.
\end{itemize}

Let us emphasize that the condition (\ref{eq:Banach-Saks}) with
$x=0$ may be rewritten for Hilbert $C^*$-modules  in the following
equivalent way
\begin{equation} \label{eq:Banach-Saks2}
    \lim_{k\rightarrow\infty}\biggl\|\frac{1}{k^2}\sum_{i,j=1}^k
    \langle x_{n(i)},x_{n(j)}\rangle\biggr\|=0 \, .
\end{equation}
Obviously, any $\tau_{\widehat{V}}$-null orthogonal sequence
$\{x_k\}$ of norm one vectors in a Hilbert $C^*$-module $V$
satisfies the condition {\rm (\ref{eq:Banach-Saks2})}.

\begin{lem}
  Condition {(mBS2)} implies condition {(mBS1)}, but
  not conversely, generally speaking.
\end{lem}

\begin{proof}
Clearly, condition {(mBS2)} implies condition {(mBS1)} as a
particular case. To see the non-equivalence of these conditions
consider the $C^*$-algebra $A=K(l_2)$ of all compact operators on
a separable Hilbert space $l_2$ as a Hilbert $C^*$-module over
itself. Then condition {(mBS1)} is fulfilled, however condition
{(mBS2)} does not hold: indeed, for any sequence of minimal
pairwise orthogonal projections $\{p_i \}$ with their least upper
bound $1_{B(l_2)}$, where $1_{B(l_2)}$ denotes the unit of the
$C^*$-algebra $B(l_2)$ of all bounded linear operators in $l_2$,
the sequence of partial sums $\left\{ q_k = \sum_{i=1}^k p_i
\right\}$ is $\tau_{\widehat{V}}$-convergent to $1_{B(l_2)}$. For
any subsequence $\{ q_{k(j)} \}$ of it the sequence $\frac{1}{j}
\sum_{i=1}^j q_{k(j)}$ converges to an infinite projection
operator on $l_2$. So the latter sequence cannot converge in norm
since infinite projection operators do not belong to $A$ which is
already norm-closed.
\end{proof}

\begin{lem}
  Any $\tau_{\widehat{V}}$-convergent sequence $\{x_n\}$ of
  a Hilbert $C^*$-module $V$ with a subsequence $\{x_{n(i)}\}$
  such that the sequence $\frac{1}{k}\sum_{i=1}^k x_{n(i)}$
  converges in norm to some element $x \in V$ admits $x$
  as its $\tau_{\widehat{V}}$-limit. Consequently, Hilbert
  $C^*$-modules with condition {(mBS2)} are sequentially
  $\tau_{\widehat{V}}$-complete.
\end{lem}

\begin{proof}
Consider any $\tau_{\widehat{V}}$-convergent sequence $\{x_n\}$ of
$V$ that admits a subsequence $\{x_{n(i)}\}$ such that the
sequence $\frac{1}{k}\sum_{i=1}^k x_{n(i)}$ converges in norm to
some $x \in V$. Then this sequence has to
$\tau_{\widehat{V}}$-converge to the same element $x$, that means
\[
      \lim_{l \to \infty} \left\| \left\langle z, \frac{1}{l}\sum_{i=1}^l
      (x_{n(i)}-x) \right\rangle \right\| \\
       =  0
\]
for any $z \in V$. Therefore, for any $z \in V$, any $\varepsilon
> 0$ there exists an $L \in  \mathbb N$ such that for any $k>l > L$
the inequality
\[
   \varepsilon >
   \left\| \left\langle z, \frac{1}{l}\sum_{i=1}^l (x_{n(i)}-x) \right\rangle -
      \left\langle z, \frac{1}{k}\sum_{i=1}^k (x_{n(i)}-x) \right\rangle \right\| =
   \frac{1}{k-l} \left\| \sum_{i=l+1}^k \langle z, (x_{n(i)}-x) \rangle \right\|
\]
holds. Selecting $k=l+1$ we arrive at
\[
   \lim_{l \to \infty} \| \langle z, (x_{n(l)}-x) \rangle \| = 0
\]
and the subsequence $\{ x_{n(l)} \}$ of the sequence $\{x_n\}$
$\tau_{\widehat{V}}$-converges to $x \in V$. Now, for any $z \in V$
we obtain
\begin{eqnarray*}
   0 & = & \lim_{n \to \infty} \lim_{l \to \infty} \| \langle z,(x_n-x_{n(l)}) \rangle \| \\
     & = & \lim_{n \to \infty} \| \langle z,(x_n-x) \rangle \|  \, .
\end{eqnarray*}
Consequently, the sequence $\{ x_n \}$ $\tau_{\widehat{V}}$-converges
to $x \in V$ by definition, and this
is equivalent to the assertion that the sequence
$\{ x_n - x \}$ $\tau_{\widehat{V}}$-converges to zero.
\end{proof}

The study of the property {(mBS1)} applied to the situation
of Hilbert $C^*$-modules gives a rather surpringly general result,
in comparison to the weak Banach-Saks property situation.

\begin{teo}\label{teo:BS1_property}
   Any Hilbert $C^*$-module has the property {(mBS1)}.
\end{teo}

\begin{proof}
We will reason in the vein of the original work \cite{B-S:30}.
Let $\{x_k\}$ be a $\tau_{\widehat{V}}$-null sequence of a
Hilbert $C^*$-module. By Proposition \ref{prop:conv_bound1}
we may suppose that the norms of all the elements $x_k$
do not exceed one. Put $n(1)=1$. By supposition for an arbitrary
$\varepsilon_1>0$ there exists a number $n(2)$ such that
$\|\langle x_{n(1)},x_i\rangle\|<\varepsilon_1$ whenever $i>n(2)$.
Analogously for an arbitrary $\varepsilon_2>0$ there is a number
$n(3)$ such that $\|\langle x_{n(2)},x_i\rangle\|<\varepsilon_2$
whenever $i>n(3)$ and so on. Set $\varepsilon_i=1/i$. Then we claim
the subsequence $\{x_{n(i)}\}$ of $\{x_k\}$ satisfies the condition
(\ref{eq:Banach-Saks2}). Indeed
\begin{align*}
  \biggl\|\frac{1}{k^2}\sum_{i,j=1}^k
    \langle x_{n(i)},x_{n(j)}\rangle\biggr\|&\le
  \frac{1}{k^2}\left(\sum_{i=1}^k \|x_{n(i)}\|^2+2\sum_{1\le i<j\le k}
    \|\langle x_{n(i)},x_{n(j)}\rangle\|\right)\\
    &\le \frac{1}{k}+\frac{2}{k^2}\left((k-1)\varepsilon_1+
    (k-2)\varepsilon_2+\dots+\varepsilon_{k-1}\right)\\
    &=\frac{1}{k}+\frac{2}{k^2}\left(k+\frac{k}{2}+\frac{k}{3}+\dots+
    \frac{k}{k-1}-(k-1)\right)\\
    &=\frac{1}{k}+\frac{2}{k^2}+\frac{2}{k}\left(\frac{1}{2}+\frac{1}{3}+\dots+
    \frac{1}{k-1}\right)=\\
    &=\frac{1}{k}+\frac{2}{k^2}+\frac{2}{k}(C-1+\ln (k-1)+o(1)),
\end{align*}
where $C\thickapprox 0,5772$ is the Euler constant and, consequently, the
right part of this estimation vanishes when $k$ goes to infinity.
\end{proof}

\begin{prop}\label{prop:Schur implies B-S}
 If a Hilbert $C^*$-module has the module Schur
 property, then it has the property {(mBS2)}.
\end{prop}

\begin{proof} Let $\{x_k\}$ be a $\tau_{\widehat{V}}$-convergent sequence
of a Hilbert $C^*$-module $V$ and let $x$ be its limit in the
completion of $V$ with respect to the topology
$\tau_{\widehat{V}}$. Then Lemma \ref{lem:seq_compl_Schur}
ensures $x$ has to belong to $V$ actually. Thus the numerical
sequence $\{ \| x_k-x \| \}$ converges to zero. Without loss of
generality (because one can pass to a subsequence) we can suppose
that $\|x_k-x\|\le 1/k^2$. This implies
\begin{align*}
    \left\|\frac{1}{n}\left(\sum_{k=1}^nx_k\right) -x\right\|&=
    \left\|\frac{1}{n}\sum_{k=1}^n(x_k -x)\right\|\\
    &\le\frac{1}{n}\sum_{k=1}^n \|x_k-x\|\\ &\le \frac{1}{n}
    \left(\sum_{k=1}^\infty \frac{1}{k^2}\right)
\end{align*}
and the right part of this expression converges to zero.
\end{proof}

\begin{teo}\label{teo:B-S for C*-alg}
  A $C^*$-algebra with a strictly positive element, in particular,
  a separable $C^*$-algebra,
  has the property {(mBS2)} if and only if it is
  unital, if and only if it has the module Schur property.
\end{teo}

\begin{proof} Theorem \ref{teo:Schur for C*-alg} and
Proposition \ref{prop:Schur implies B-S} ensure that a unital
$C^*$-algebra has the property {(mBS2)}. To check
the converse implication let us fix an arbitrary non-unital
$C^*$-algebra $A$ with a strictly positive element and consider
$V=A$ as a Hilbert $A$-module. Then exactly by the same reasons
as explained in the proof of Theorem \ref{teo:Schur for C*-alg}
any element of the right multiplier algebra $RM(A)$ may be obtained
just as a limit  of a sequence from $A$ with respect to the right
strict topology, which coincides with the topology
$\tau_{\widehat{V}}$. Therefore, one can find a
sequence $\{x_k\}$ of $A$ which converges with respect to the
topology $\tau_{\widehat{V}}$ (for instance, to the unit of
$RM(A)$), but which does not converge in norm. Then for any
subsequence $\{x_{k(i)}\}$ of $\{x_k\}$ the sequence
$\{\frac{1}{n}\sum_{i=1}^n x_{k(i)}\}$ of $A$ has obviously the
same $\tau_{\widehat{V}}$-limit in $RM(A)$ as the sequence
$\{x_k\}$, so $\{\frac{1}{n}\sum_{i=1}^n x_{k(i)}\}$ cannot
converge in norm. Thus $A$ does not have the property {(mBS2)}.
\end{proof}

%\begin{cor}
%   Any separable unital $C^*$-algebra has the
%   property {(mBS2)}.
%\end{cor}

\begin{teo}
  Assume a Hilbert $C^*$-module $V$ admits a decomposition into
  a direct sum $V=V_1\oplus V_2$. Then the following conditions
  are equivalent:
  \begin{enumerate}
    \item $V$ has the property {(mBS2)}.
    \item Both $V_1$ and $V_2$ have the property {(mBS2)}.
  \end{enumerate}
\end{teo}

\begin{proof}
Obviously, (i) implies (ii) and we have just to ensure the inverse
conclusion. Consider any $\tau_{\widehat{V}}$-convergent sequence
$\{x_i\}$ of $V$. Then its elements admit decompositions as
$x_i=x_i^{(1)}\oplus x_i^{(2)}$ with $x_i^{(1)}\in V_1$,
$x_i^{(2)}\in V_2$. Both the sequences $\{x_i^{(j)}\}$ are
$\tau_{\widehat{V_j}}$-convergent sequences  ($j=1, 2$).
Therefore, there are subsequences $\{x_{i(k)}^{(j)}\}$ of
$\{x_i^{(j)}\}$ such that their means $\frac{1}{n}\sum_{k=1}^n
x_{i(k)}^{(j)}$ converge in norm to some elements $x^{(j)}$ of
$V_j$, respectively for $j=1,2$. Thus, for the sequence $\{
x_{i(k)} \}$ and for $x=x^{(1)}\oplus x^{(2)}$ we have the
following estimate:
\begin{align*}
   \biggl\|\frac{1}{n}\left(\sum_{k=1}^n x_{i(k)}\right)-(x^{(1)}\oplus x^{(2)})\biggr\|
   &=\biggl\|\left(\frac{1}{n}\sum_{k=1}^n x_{i(k)}^{(1)}-x^{(1)}\right)\oplus
   \left(\frac{1}{n}\sum_{k=1}^n x_{i(k)}^{(2)}-x^{(2)}\right)\biggr\|\\
   &\le \biggl\|\frac{1}{n}\sum_{k=1}^n x_{i(k)}^{(1)}-x^{(1)}\biggr\|+
 \biggl\|\frac{1}{n}\sum_{k=1}^n x_{i(k)}^{(2)}-x^{(2)}\biggr\|
\end{align*}
Since the upper bounds converge to zero as $k$ tends to infinity, the
sum in the first term converges to $x$ in norm. So $V$ has the
property {(mBS2)}.
\end{proof}

\begin{cor}
  Let $A$ be a $\sigma$-unital $C^*$-algebra
  with the property {(mBS2)}.
  Then finitely generated Hilbert modules over
  unital $A$ and finitely generated projective
  Hilbert modules over $A$ have the property
  {(mBS2)}.
\end{cor}

Let $A$ be a $C^*$-algebra and consider the infinite
matrix algebra $M_\infty(A) =\cup_{n=1}^\infty M_n(A)$.
For vectors $x_1,\dots,x_n$ of a Hilbert $A$-module $V$
consider the Gram matrix based on them:
\[
  G(x_1,\dots,x_n)=
  \begin{pmatrix}
  \langle x_1,x_1\rangle & \langle x_1,x_2\rangle\dots & \langle x_1,x_n\rangle \\
  \langle x_2,x_1\rangle & \langle x_2,x_2\rangle\dots & \langle x_2,x_n\rangle \\
  \hdotsfor{3} \\
  \langle x_n,x_1\rangle & \langle x_n,x_2\rangle\dots & \langle
  x_n,x_n\rangle
  \end{pmatrix} \, ,
\]
with entries $g_{ij}=\langle x_i,x_j\rangle$, besides
$G(x_1,\dots,x_n)\in M_\infty(A)$. Then the condition
(\ref{eq:Banach-Saks2}) on a sequence $\{x_i\}$ of $V$
with $x=0$ can be equivalently reformulated as the
assertion that $\{\frac{1}{k^2}\sum_{i,j=1}^k g_{ij}\}$
forms a sequence that converges to zero with respect
to the norm topology. Instead of condition
(\ref{eq:Banach-Saks2}) one can require that the sequence
$\{\frac{1}{n^2} G(x_1,\dots,x_n)\}$ of $M_\infty(A)$
has to converge to zero in norm. For instance, this
condition holds whenever $\{x_i\}$ is an orthogonal system
of norm one vectors, because $\|G(x_1,\dots,x_n)\|\le
\sum_{i,j=1}^n \|g_{ij}\|$ (cf.~\cite{Murphy}). Let us
remark that in general neither the value $\|\frac{1}{n^2}
\sum_{i,j=1}^n \langle x_{i},x_{j}\rangle\|$ nor the value
$\|\frac{1}{n^2} G(x_1,\dots,x_n)\|$ majorize each other.

%%%%%%%%%%%%%%%%%%%%%%%%%%%%%%%%%%%%%%%%%%%%%%%%%%%%%%%%
\section{Dual modules and module Schur and Banach-Saks properties}
%%%%%%%%%%%%%%%%%%%%%%%%%%%%%%%%%%%%%%%%%%%%%%%%%%%%%%%%
In this section we establish an interrelation between
self-duality of countably generated Hilbert $C^*$-modules
and their properties to have the module Schur or the
module Banach-Saks property, respectively. The main
goal is an alternative characterization of $C^*$-dual
Hilbert $C^*$-modules of representatives of this
class.

Suppose $X$ denotes the completion of $V$ with respect to the
topology $\tau_{\widehat{V}}$, and $X_s$ denotes the subset of $X$
consisting of all equivalence classes of
$\tau_{\widehat{V}}$-Cauchy sequences from $V$. (Two Cauchy
sequences are supposed to be equivalent if and only if their
difference is a $\tau_{\widehat{V}}$-null sequence.) Let us show
how the set $X_s$ can be canonically identified with a subset of
the dual module $V'$, similarly like $V$ is canonically embeddable
into its $C^*$-dual $A$-module $V'$. Let $\{x_i\}$ be a
$\tau_{\widehat{V}}$-Cauchy sequence from $V$, i.e. the sequence
$\{\langle x_i, z\rangle\}$ converges in norm to some element
$x(z)\in A$ for any $z\in V$ since $\|\langle u,v \rangle \|
= \|\langle v,u \rangle \|$ for any $u,v \in V$. As we see this
construction defines an $A$-linear functional $\langle x, \cdot
\rangle : V\rightarrow A$. Now, for this functional $\langle x,
\cdot \rangle$, for $z\in V$ and for any $\varepsilon>0$ there
exists a number $i_0$ such that
\[
  \|x(z)\|\le \|\langle x_i,z\rangle\|+\varepsilon\le
  \|x_i\|\|z\|+\varepsilon
\]
whenever $i>i_0$, where $\|x_i\|<C$ for some constant $C$ and
for all $i$ by Proposition \ref{prop:conv_bound1}. Thus, the
functional $\langle x, \cdot \rangle$ is
bounded and belongs to the dual module $V'$.

\begin{prop}
   The topology $\tau_{\widehat{V}}$ has a countable base of
   neighborhoods of zero provided
   $V$ is a countably generated Hilbert $A$-module.
\end{prop}

\begin{proof} Let $V=\overline{\mathrm{span}}_A \{x_i : i\in
\mathbb{N}\}$, i.e. the sequence $\{x_i\}_{i=1}^\infty$ generates
$V$ over $A$. Then we claim the countable family of sets
\begin{gather}\label{eq:count_base}
   U_{x_{i_1},\dots,x_{i_n},1/N}=\{y\in V : \|\langle
   x_{i_1},y\rangle\|<1/N,\dots,\|\langle x_{i_n},y\rangle\|<1/N\},\\
   1<i_1,\dots,i_n<\infty, n,N \in \mathbb{N}\notag
\end{gather}
forms the base of neighborhoods of zero for the topology
$\tau_{\widehat{V}}$. Indeed, let us take into consideration any
$\tau_{\widehat{V}}$-neighborhood of zero of the form
\[
U_{x,\varepsilon}=\{y\in V : \|\langle x,y\rangle\|<\varepsilon\},
\]
where $x\in V$, $\varepsilon>0$. It is just enough to check that
any of these open sets contains one of the sets (\ref{eq:count_base}).
By supposition, for any $\delta>0$ there are $a_1,\dots,a_n$ of $A$
such that
\begin{gather*}
    \biggl\|x-\sum_{i=1}^n x_ia_i\biggr\|<\delta.
\end{gather*}
Consider a neighborhood $U_{x_{1},\dots,x_{n},1/N}$. By
Proposition \ref{prop:conv_bound3} we can assume that the norms of
all the elements of the generating set do not exceed one, whenever
$N$ is large enough.
Therefore for any $z\in U_{x_{1},\dots,x_{n},1/N}$ one can deduce
\begin{align*}
    \|\langle x,z\rangle\|&\le \biggl\|\biggl\langle x-\sum_{i=1}^n
    x_ia_i,z\biggr\rangle\biggr\|+ \biggl\|\biggl\langle\sum_{i=1}^n
    x_ia_i,z\biggr\rangle\biggr\|\\
    &\le \delta\|z\|+\sum_{i=1}^n\|a_i\|\|\langle x_i,z\rangle\|\\
    &\le \left( \delta+\frac{1}{N}\sum_{i=1}^n\|a_i\| \right) \|z\| \, .
\end{align*}
The numerical expression turns out to be less than the selected
$\varepsilon$, whenever both $\delta$ and $1/N$ are small enough
and $n \in \mathbb N$ is fixed. So we are done.
\end{proof}

\begin{cor}
   Let $V$ be a countably generated Hilbert $C^*$-module. Then
   under the notations above $X_s$ coincides with the completion $X$
   of $V$ with respect to the topology $\tau_{\widehat{V}}$ and,
   consequently, $X$ can be isometrically embedded into the
   the $C^*$-dual module $V'$ extending the canonical embedding of
   $V$ into $V'$.
\end{cor}

\begin{teo}\label{teo:X=V'_count_generated}
   Let $V$ be a countably generated Hilbert module over a unital
   $C^*$-algebra $A$.
   Then the completion $X$ of $V$ with respect to the topology
   $\tau_{\widehat{V}}$ coincides with the dual module $V'$
   extending the canonical embedding of $V$ into $V'$.
\end{teo}

\begin{proof} We already know that $X$ is included into $V'$ via
the canonical embedding of $V$ into $V'$ and should just check the
coincidence in the particular case. Firstly, let us suppose $V$ is
the standard module $l_2(A)$. Then we can use the characterization
of the dual module $l_2(A)'$ described in Theorem \ref{prop:dual
for the standard module}. For any $\beta=(b_i)\in l_2(A)'$, where
$\|\sum_{i=1}^n b_i^*b_i\|\le C$ for all $n$ and some finite
constant $C$, let us consider the finite vectors
$\alpha_n=(b_1,\dots,b_n,0,0,\dots) \in l_2(A)$. Then for any
$y=(y_i)\in l_2(A)$ we have
\begin{align*}
    \|\langle \alpha_n-\alpha_m, y\rangle\|&=\biggl\|\sum_{i=n+1}^m
    b_i^*y_i\biggr\|\\&\le \biggl\|\sum_{i=n+1}^m
    b_i^*b_i\biggr\|^{1/2}\biggl\|\sum_{i=n+1}^m
    y_i^*y_i\biggr\|^{1/2}\\
    &\le C^{1/2}\biggl\|\sum_{i=n+1}^m y_i^*y_i\biggr\|^{1/2}
\end{align*}
and this expression goes to zero provided $m$, $n$ go to infinity.
Thus, $\{\alpha_n\}$ is a $\tau_{\widehat{V}}$-Cauchy sequence
and, moreover, its $\tau_{\widehat{V}}$-limit, obviously,
coincides with $\beta$. This shows the desirable identification of
$X$ and $V'$ that extends the canonical embedding of $V$ into
$V'$.

Now assume $V$ is an arbitrary countably generated Hilbert
$A$-module. Then by the Kasparov's stabilization theorem there
exists a Hilbert $A$-module $W$ such that $l_2(A)=V\oplus W$. Let
us denote by $P_V : l_2(A)\rightarrow l_2(A)$ the corresponding
orthogonal projection onto the first summand, so
$\mathrm{Range}(P_V)=V$. Now for any functional $f\in V'$ one can
define a functional $\widetilde{f}\in l_2(A)'$ as the composition
$\widetilde{f}=fP_V$. Because $l_2(A)'$ coincides with the
$\tau_{\widehat{l_2(A)}}$-completion of $l_2(A)$ there is a
sequence $\{z_i\}$ of $l_2(A)$, where $z_i=x_i\oplus y_i$, such
that its $\tau_{\widehat{l_2(A)}}$-limit is $\widetilde{f}$.
Thus, the sequence $\{x_i\}$ of $V$ converges to $f$ with respect
to the topology $\tau_{\widehat{V}}$.
\end{proof}

\begin{rk}\rm
Theorem \ref{teo:X=V'_count_generated} is true
also for $\sigma$-unital $C^*$-algebras $A$, but the proof
involves more complicated techniques and is left to the reader.
Actually, it is the corollary of the crucial
general result \cite[Theorem 6.4]{FrankPOS} and Proposition
\ref{prop:conv_bound1}. Let us also stress for $V=A$ that in the
non-unital, $\sigma$-unital case the theorem above is equivalent
to the well-known fact that the completion of $A$ with respect
to the left strict topology coincides with the left multiplier
algebra $LM(A)$.
\end{rk}

Applying Theorem \ref{teo:X=V'_count_generated} we get:

\begin{cor}
  Let $V$ be a countably generated Hilbert module
  over a unital $C^*$-algebra.
  Then $V$ is self-dual whenever it has the module
  Schur property.
\end{cor}

\begin{cor}\label{cor:equivalence_BS_properties}
  Let $V$ be a countably generated self-dual Hilbert
  module over a unital $C^*$-algebra.
  Then $V$ has the property {(mBS1)} if and only if
  it has the property {(mBS2)}.
\end{cor}

\begin{proof} Suppose the property {(mBS1)} holds for $V$.
Consider any $\tau_{\widehat{V}}$-convergent sequence $\{x_i\}$ of
$V$. Then by supposition and by Theorem \ref{teo:X=V'_count_generated}
its $\tau_{\widehat{V}}$-limit $x$ belongs to $V$. Then the
sequence $\{y_i=x_i-x\}$ is a $\tau_{\widehat{V}}$-null sequence.
Hence, it admits a subsequence $\{y_{i(k)}\}$ satisfying the
equality (\ref{eq:Banach-Saks}), so the sequence
$\frac{1}{n}\sum_{k=1}^n x_{i(k)}$ converges to $x$ in norm.
\end{proof}

\begin{cor}
  Any countably generated self-dual Hilbert module
  over a unital $C^*$-algebra possesses the property {(mBS2)}.
\end{cor}

The proof follows from Theorem \ref{teo:BS1_property} and
Corollary \ref{cor:equivalence_BS_properties}.

As a final result we can state, that the module Schur and
the {(mBS2)} properties are different in general,
because any standard Hilbert module over a finite dimensional
$C^*$-algebra is self-dual and, consequently, has the
property {(mBS2)}, however it does not have the module
Schur property by Proposition \ref{prop:Schur stand mod}.

\medskip \noindent
{\bf Acknowledgement:} The presented work is part of the research
project ''$K$-Theory, $C^*$-Algebras, and Index Theory'' of
Deutsche Forschungsgemeinschaft (DFG). The authors are grateful to
DFG for the support. The essential part of the research was done
during the visit of the second author to the Leipzig University of
Applied Sciences (HTWK), and he appreciates its hospitality a lot.

The authors are grateful to A. Ya. Helemskii for
helpful discussions.

%\bibliographystyle{plain}
%\bibliography{references08}
%%%%%%%%%%%%%%%%%%%%%%%%%%%%%%%%%%%%%%%%%%%%%%%%%%%%%%%%%%%%%%%%%%%%%%%
%%%%%%%%%%%%%%%%%%%%%%%%%%%%%%%%%%%%%%%%%%%%%%%%%%%%%%%%%%%%%%%%%%%%
%%%%%          The bibliography
%%%%%%%%%%%%%%%%%%%%%%%%%%%%%%%%%%%%%%%%%%%%%%%%%%%%%%%%%%%%%%%%%%%%%

\end{document}